\documentclass[12pt]{article}
\begin {document}

\newenvironment{proof}{
\par
 {\large Proof.}\rm}

\title{Simple decompositions of the exceptional Jordan algebra }
\author{
M.V.Tvalavadze \thanks{supported by NSERC Grant 227060-00}\\
 }
\date{}
\maketitle \sloppy

\noindent{Abstract.} {\footnotesize This paper presents some
results on simple exceptional Jordan algebra over algebraically
closed field $\Phi$ with characteristic not 2. Namely an example
of simple decomposition of $H(O_3)$ into the sum of two
subalgebras of the type $H(Q_3)$ is produced, and it is shown that
this decomposition is the only possible in terms of the type of
simple subalgebras.}

\par\medskip
\par\medskip

{\bf 1. Introduction.} Let $\Phi$ be an algebraically closed field
with characteristic unequal to 2. Then all composition algebras
are split in a sense of having no zero divisors. Therefore, they
are uniquely determined by their dimensions. In particular, there
exists only one octonion algebra over $\Phi$ denoted as $O$. The
algebra of split octonions can be represented as a direct vector
sum of $Q$ and $vQ$ where $Q$ is the split quaternion algebra and
$v^2=1$, $O=Q\oplus Qv$. Multiplication in $O$ is defined by
$$ (q+rv)(s+tv)=(qs+\bar tr)+(tq+r\bar s)v$$
for $q,r,s,t$ in $Q$. The involution in $O$ is
$$x=q+rv\to \bar x=\bar q-rv.$$
Therefore, $O$ is a vector space of dimension 8 and degree 2. Next
we consider the set $H_3(O)$ of Hermitian $3\times 3$ matrices
over $\Phi$. An arbitrary element $A$ from this set has the form
$$ A=\left(\begin{array}{ccc}
          \alpha  &   x    &  y \\
          \bar x  & \beta  &  z \\
          \bar y  & \bar z &\gamma
          \end{array}\right),
$$
where $\alpha$, $\beta$, $\gamma \in \Phi$, $x,y,z \in O$.
Obviously, $H_3(O)$ is a vector space of dimension 27. If $x[ij]$
denotes $xE_{ij}+\bar x E_{ji}$ where $x\in O$, $E_{ij}$ a
standard matrix unit, then $A=\alpha E_{11}+\beta E_{22}+\gamma
E_{33}+x[12]+y[13]+z[23]$.

Let $O$ be an octonion algebra, " --- " the canonical involution.
Let $a,b,c,d,x,z\in O$; recall \cite{jac}

(1) $x^2-t(x)x+n(x)=0$ where

(2) $n(x)=x\bar x=\bar x x \in \Phi$

(3) $t(x)=t(x,1)$ where

(4) $t(a,b)=n(a+b)-n(a)-n(b)=a\bar b+b\bar a=\bar a b+\bar b
a=t(\bar a,b)$

(5) $\bar a (ab)=n(a)b=(ba)\bar a$.

In \cite{blij} Van Der Blij and Springer obtained the standard
basis $\{x_i,y_i| 0\le i\le 3\}$ of the split octonions which
satisfies the following identities:
$$ x_0x_i=x_i,\quad y_0y_i=y_i,\quad 0\le i\le 3$$
$$ x_iy_i=-x_0,\quad x_ix_{i+1}=y_{i+2},\quad y_iy_{i+1}=x_{i+2},\quad 1\le i\le 3$$
where the indices are taken modulo 3 and all other products are
zero or obtained by applying  the canonical involution to one of
the above.
\par\medskip

{\bf 2. Example of a simple decomposition of $H(O_3)$.} To exhibit
an example of a simple decomposition of $H(O_3)$ we need the
following Lemma.

{\large Lemma 1.} {\it The set $\cal S$ consisting of all
hermitian matrices of the form
$$     \left(\begin{array}{ccc}
          \alpha  &   xv    &  yv \\
          -xv  & \beta  &  z \\
          -yv  & \bar z &\gamma
          \end{array}\right), \eqno (1)
$$
where $x,y,z\in Q$ is a simple algebra of the type $H(Q_3)$}.
\par\medskip
\begin{proof} First, it is easy to verify that $\cal S$ is closed
with respect to the product $X\circ Y=\frac{1}{2}(XY+YX)$ where
$XY$ is the ordinary matrix product in $O_3$. Therefore, $\cal S$
is a Jordan subalgebra of $H(O_3)$.

Next we are going to show that $\cal S$ has no proper non-zero
ideals. Assume the contrary, that is, there exists an ideal $I$
such that $I\ne\{0\}$, $I\ne {\cal S}$. Hence we can select a
non-zero $a_0$ from $I$ of the form:
$a_0=\alpha_0E_{11}+\beta_0E_{22}+\gamma_0E_{33}+x_0v[12]+y_0v[13]+z_0[23]$.
Let one of off-diagonal elements be non-zero. For clarity, $x_0\ne
0$. Then, multiplying $a_0$ by $E_{11}$ and $E_{22}$ (in the
Jordan sense) we obtain $a_1=x_0v[12]$ which is also an element of
$I$. First we assume that  $n(x_0)=x_0\bar x_0=0$. Then for any
$b=\alpha E_{11}+\beta E_{22}+\gamma E_{33}+x[12]+y[13]+z[23]$ and
$c=\alpha' E_{11}+\beta' E_{22}+\gamma'
E_{33}+x'[12]+y'[13]+z'[23]$
$$b\circ c=\alpha''E_{11}+\beta''E_{22}+\gamma''E_{33}+b_{12}[12]+b_{13}[13]+b_{23}[23] \eqno (2)$$
where
$$b_{12}=(\alpha x'+\beta' x+y z'+\alpha' x+\beta  x'
+y' z)v,$$
$$b_{13}=(\alpha y'+x \bar z'+\gamma' y+\alpha' y+x'\bar
z+\gamma y')v,$$
$$b_{23}=\bar y' x+\beta z'+\gamma' z+\bar y x' +\beta' z+z'\gamma,$$
$$\alpha''=2\alpha\alpha'+t(x,x')+t(y,y'),$$
$$\beta''=2\beta\beta'+t(x,x')+t(z,z'),$$
$$\gamma''=2\gamma\gamma'+t(z',z)+t(y,y').$$

 Set $c=a_1=x_0v[12]$. If we choose $b$ such that $\alpha=0$ and the other
 coefficients of $b$
are non-zero, then $b\circ c$ has the form (2) where $b_{12}=\beta
x_0v$, $b_{13}=x_0\bar z v$, $b_{23}=\bar y x_0$,
$\alpha''=\beta''=t(x_0,x)$, $\gamma''=0$. Again, multiplying
$b\circ c$ by $d=x''[12]+ z''[13]$, $x'',z''\ne 0$  we obtain
$(b\circ c)\circ
d=\alpha'''E_{11}+\beta'''E_{22}+\gamma'''E_{33}+c_{12}[12]+c_{13}[13]+c_{23}[23]$
such that $c_{12}=(\alpha''
x''+b_{13}z''+\beta''x'')v=(2t(x_0,x'')x''+x_0\bar z z'')v$. Next
we are going to show that we can choose $x''$ in such a way that
$n(2t(x_0,x'')x''+x_0\bar z z'')\ne 0$. Since $z''$ is an
arbitrary element, we can set $z''=0$. Consider
$r=2t(x_0,x'')x''$. Since $t(x_0,x'')$ is a non-singular bilinear
form and $x_0\ne 0$, there exists $q\in Q$: $t(x_0,q)\ne 0$. If
$n(q)\ne 0$, then $n(2t(x_0,q)q)=4t(x_0,q)^2n(q)\ne 0$. Therefore,
    $x''=q$ is a required element. If $n(q)=0$, then we consider
$q+\delta$ where $\delta\in \Phi$. Notice that
$n(q+\delta)=n(q)+\delta t(q)+\delta^2=\delta(t(q)+\delta)\ne 0$
if $\delta\ne 0,\, -t(q)$, and $t(x_0,q+\delta)=t(x_0,q)+\delta
t(x_0)$. If $t(x_0)=0$, then $t(x_0,q+\delta)\ne 0$ for any
$\delta \in \Phi$. If $t(x_0)\ne 0$, then $t(x_0,q+\delta)\ne 0$
for any $\delta\ne t(x_0,q)/t(x_0)$. Hence, if $\delta \ne 0,\,
-t(q),\, t(x_0,q)/t(x_0)$, then $n(2t(x_0,q+\delta), q+\delta)\ne
0 $. Therefore, $x''=q+\delta$ is a required element. As shown
above we can always choose $x''$ such that $n(c_{12})\ne 0$. Then,
performing multiplication of $(b\circ c)\circ d$ by $E_{11}$ and
$E_{22}$ we come to $\bar a_0=c_{12}v[12]$ where $n(c_{12})\ne 0$.

These considerations allow us to assume that $n(x_0)\ne 0$. Then
$a_0^2=-n(x_0)(E_{11}+E_{22})$, $n(x_0)\in\Phi$, $n(x_0)\ne 0$.
Hence $E_{11}+E_{22}\in I$. It follows that for any $\alpha,\beta
\in\Phi$, $x\in Q$: $\alpha E_{11}+\beta E_{22}+xv[12]\in I$.
Then, it is easily seen that $I={\cal S}$, a contradiction.

As a result, $\cal A$ is a simple special Jordan algebra of degree
3 and dimension 15. Hence, ${\cal S}\cong H({\cal Q}_3)$. This
proves the Lemma.
\end{proof}

Before we continue our discussion we need a few facts concerning
$\mbox{Aut}\, H(O_3)$. Let $\{ e_1, e_2, e_3\}$ be a reducing set
of idempotents of $H(O_3)$, and let $H(O_3)=\sum_{ij} J_{ij}$ be
the corresponding Pierce decomposition. Set ${\cal C}=\Phi
(e_1-e_2)\oplus J_{12}$, $\dim {\cal C}=9$. Then $\cal C$ is an
algebra of the bilinear form $Q(x,y)$ where
$Q(x)=\frac{1}{2}t(x^2)$. From \cite{jac}, if $v_1, \ldots, v_{2r}
\in {\cal C}$, $\Pi_i Q(v_i)=1$, then
$\eta=U_{\theta(v_1)}U_{\theta(v_2)}\ldots U_{\theta(v_{2r})}$
where $\theta(x)=e_3+x$ is an automorphism of $H(O_3)$ that leaves
$e_3$ fixed.

Let $e_3=E_{33}$, and
$$a=\left(\begin{array}{ccc}
           1/\sqrt 3    & (1+iv)/\sqrt 3 &   0 \\
          (1-iv)/\sqrt 3  & -1/\sqrt 3  &    0 \\
           0  & 0 & 0
          \end{array}\right),\eqno (3)
$$
$$b=\left(\begin{array}{ccc}
      1& 0& 0\\
      0& -1& 0\\
      0& 0& 0
      \end{array}\right). \eqno (4)
$$
Obviously, $a,b\in {\cal C}$, and  $Q(a)=1$, $Q(b)=1$, therefore
$\varphi=U_{\theta(a)}U_{\theta(b)}$ is an automorphism of
$H(O_3)$.
\par\medskip
{\large Example 1.} {\it There exists a simple decomposition of
$H(O_3)$ of the form $H(O_3)={\cal A}+{\cal B}$ where $\cal A$,
${\cal B}\cong H(Q_3)$. }
\par\medskip
\begin{proof}
First we notice that by \cite{jac} any two subalgebras of the type
$H(Q_3)$ are conjugate under automorphism of $H(O_3)$. Therefore,
we can assume that the first subalgebra has the canonical form
$${\cal A}=\left\{ \left(\begin{array}{ccc}
          \alpha  &   x    &  y \\
          \bar x  & \beta  &  z \\
          \bar y  & \bar z &\gamma
          \end{array}\right)\right\},
$$
where $\alpha$, $\beta$, $\gamma \in \Phi$, $x,y,z \in Q$. Next we
choose $a$ of the form (3) and $b$ of the form (4). As shown
above, $\varphi=U_{\theta(a)}U_{\theta(b)}$ is an automorphism of
$H(O_3)$. Recall $xU_{\theta(y)}=2(x\circ \theta(y))\circ
\theta(y)-x\circ \theta(y)^2$. Notice
$(a+E_{33})^2=(b+E_{33})^2=E_{11}+E_{22}+E_{33}$. Let ${\cal
B}=\varphi({\cal S})$ where $\cal S$ is an algebra of the form
(1). Any $b\in {\cal B}$ has the form $b=\beta_1 E_{11}+\beta_2
E_{22}+\beta_3 E_{33}+\sum_{ij}\bar b_{ij}[ij]$ where
$$\bar b_{12}=xv-\alpha (1-iv), \quad\alpha\in\Phi$$
$$\bar b_{13}=1/\sqrt 3 \bar z-((1+\frac{2}{\sqrt 3}) y +1/\sqrt 3\bar z i)v$$
$$\bar b_{23}=-(1-\frac{2}{\sqrt 3})\bar z+\frac{2}{\sqrt
3}yv(1+iv).$$ It is easy to verify that $H(O_3)={\cal A}+{\cal
B}$.
\end{proof}
{\large Remark}. The exceptional Jordan algebra can be represented
as the sum of two special algebras.

In \cite{R1} Racine has determined the maximal (unital)
subalgebras of finite-dimensional special simple linear Jordan
algebras. The paper \cite{R2} completes the determination of the
maximal subalgebras of finite-dimensional simple Jordan algebras
by considering the exceptional algebras. The main result in
\cite{R2} is the following classification Theorem.

{\large Theorem 1.}(Racine)  {\it A subalgebra ${\cal B}$ of a
finite-dimensional simple exceptional central quadratic Jordan
algebra ${\cal J}$ is maximal if and only if it is isomorphic to
an algebra of the following types:

(1) ${\cal A}^{(+)}$, ${\cal A}$ an associative division algebra
of degree 3 over its center $\Phi$.

(2) $H({\cal A}, *)$, ${\cal A}$ an associative division algebra
of degree 3 over its center $E$, an involution of the second kind
with $\Phi=H(E, *)$.

(3) $H({\cal L}_3, J)$, ${\cal L}$ is a division quaternion
algebra.

(4) ${\cal J}_1(e)\oplus {\cal J}_0(e)=\Phi e\oplus {\cal
J}_0(e)$, $e$ a primitive idempotent.

(5) $H({\cal S}_3)$, ${\cal S}=\Phi x_0\oplus \Phi x_1\oplus \Phi
x_2\oplus\Phi y_0\oplus\Phi y_1\oplus\Phi y_3$.

(6) ${\cal J}(\Phi z)$, the idealizer of $\Phi z$. }
\par\medskip
 Note that
their dimensions are 9, 9, 15, 11, 21 and 18, respectively. The
first three are special simple, while the forth is semisimple and
special. The last two are exceptional. Notice the Wedderburn
factor of $H({\cal S}_3)$ is $H({\cal L}_3)$ where $\cal L$ is a
quaternion algebra, and the Wedderburn factor of ${\cal J}(\Phi
z)$ is $S=\Phi e_1\oplus \Phi e_2\oplus (\Phi x_1\oplus \Phi
x_2\oplus \Phi x_3\oplus \Phi y_1\oplus \Phi y_2\oplus \Phi
y_3)[12]$ (we use here the basis of \cite{blij}).

{\large Theorem 2.} {\it Decomposition in Example 1 is the  only
possible simple decomposition in $H(O_3)$ in terms of the types of
simple subalgebras.}
\par\medskip
\begin{proof}
Assume the contrary, that is,  there exists a simple decomposition
different from one in Example 1. Let $H(O_3)={\cal A}+{\cal B}$
where $\cal A$, $\cal B$ simple subalgebras, at least one of them
not of the type $H(Q_3)$. Recall that $\dim\,H(F_3)=6$,
$\dim\,F_3^{(+)}=9$, $\dim\,H(Q_3)=15$. Then, by the dimension
arguments, one of the above subalgebras, for example $\cal B$,
should have the type $J(f,1)$ ( the algebra of a bilinear  form
$f$ with the identity 1).

According to the Racine's classification of maximal subalgebras,
$\cal B$ cannot be  maximal in $H(O_3)$. Therefore, we can cover
$\cal B$ with some maximal subalgebra $\cal M$ of one of the above
types. If $\cal M$ has type (1) or (2),  then $\dim\,{\cal  M}=9$.
Hence, $\dim\, {\cal B}\le 9$. Likewise if $\cal M$ has the type
(4), then $\dim\,{\cal M}\le 11$. Therefore, $\dim\, {\cal B}\le
11$.

Now if ${\cal M}\cong H(Q_3)\cong H(F_6, j)$, then $\dim\,{\cal
B}\le 4$ (see \cite{m}).

For the last two possibilities for $\cal M$ we recall that
$H({\cal S}_3)$ and ${\cal J}(Fz)$ are both exceptional Jordan
algebras with the non-zero radicals:
$$H({\cal S}_3)=H(Q_3)\oplus R,$$
$${\cal J}(Fz)=S\oplus R $$
where $S=\Phi e_1\oplus \Phi e_2\oplus (\Phi x_1\oplus \Phi
x_2\oplus \Phi x_3\oplus \Phi y_1\oplus \Phi y_2\oplus \Phi
y_3)[12]$, $\dim\, S\le 8$. By Malcev Theorem (see \cite{Mc2}) if
${\cal B}\subseteq {\cal M}=S\oplus R$, there exists an
automorphism $\varphi: {\cal M}\to {\cal M}$, such that
$\varphi({\cal B})\subseteq S$. Hence, in the case of ${\cal
M}=H({\cal S}_3)$, $\dim\,{\cal B}\le 4$, and in the case of
${\cal M}={\cal J}(\Phi z)$, $\dim\,{\cal B}\le 8$.

As a result if $\cal B$ of the type $J(f,1)$ is a subalgebra of
$H(O_3)$, then $\dim\,{\cal B}\le 11$. Further, $\cal A$ has
either the type $J(f,1)$ or $H(Q_3)$. Hence, $\dim {\cal A}\le
15$. Thus, by the dimension arguments, the decomposition does not
hold. Theorem is proved.
\end{proof}

\par\medskip
{\large Acknowledgement.} The author uses this opportunity to
thank her supervisor Prof. Bahturin for his helpful cooperation,
many useful  ideas and suggestions.

\par\medskip
\par\medskip
{\it Department of Mathematics and Statistics

Memorial University of Newfoundland

St. John's, NL

Canada,

e-mail: marina@math.mun.ca }

\end{document}